\documentclass[10pt]{amsart}
\usepackage{amssymb}
\numberwithin{equation}{section}

\newtheorem{theorem}{Theorem}

\newtheorem{proposition}[theorem]{Proposition}
\newtheorem{lemma}[theorem]{Lemma}
\theoremstyle{definition}

\theoremstyle{remark}
\newtheorem{remark}{Remark}
\newtheorem{example}{Example}

\newcommand{\R}{\mathbb{R}}
\newcommand{\C}{\mathbb{C}}
\newcommand{\Z}{\mathbb{Z}}

\newcommand\lie[1]{\mathfrak{#1}}

\newcommand{\fh}{\lie{h}}

\newcommand{\ft}{\lie{t}}

\def	\inv	{^{-1}}

\newcommand\Diff{\mathop{\it Diff}\nolimits}
\newcommand\Pol{\mathop{\it Pol }\nolimits}
\begin{document}

\title{On maximal tori in the contactomorphism groups of regular contact 
manifolds}

\author{Eugene Lerman}
\address{Department of
Mathematics, University of Illinois, Urbana, IL 61801}
\email{lerman@math.uiuc.edu}

\thanks{Supported by the Swiss NSF, US NSF grant DMS-980305 and by 
R. Kantorovitz.}

\begin{abstract}
By a theorem of Banyaga the group of diffeomorphisms of a manifold $P$
preserving a regular contact form $\alpha$ is a central $S^1$
extension of the commutator of the group of symplectomorphisms of the
base $B = P/S^1$.  We show that if $T$ is a Hamiltonian maximal torus
in the group of symplectomorphism of $B$, then its preimage under the
extension map is a maximal torus not only in the group $\Diff(P,
\alpha)$ of diffeomorphisms of $P$ preserving $\alpha$ but also in the
much bigger group of contactomorphisms $\Diff (P, \xi)$, the group of
diffeomorphism of $P$ preserving the contact distribution $\xi = \ker
\alpha$.  We use this (and the work of Hausmann, and Tolman on polygon
spaces) to give examples of contact manifolds $(P, \xi = \ker \alpha)$
with maximal tori of different dimensions in their group of
contactomorphisms.
\end{abstract}

\maketitle
\section{Introduction and the main result}

Let $(B, \omega)$ be an integral symplectic manifold.  Assume also
that $B$ is compact and simply connected.  Since the class
$[\omega]\in H^2 (B, \R)$ is integral, there exists a principal circle
bundle $S^1 \to P \stackrel{\pi}{\to} B$ with Euler class $[\omega]$.
Moreover there exists a connection 1-form $\alpha$ on $P$ with
$d\alpha = \pi^* \omega$ \cite{BW}. Consequently $\alpha$ is a contact
form and $\xi := \ker \alpha$ is a contact distribution.  In this note
we exploit a relationship between the group $\Diff(P, \xi)$ of
diffeomorphisms of $P$ preserving the contact distribution and the
group of symplectomorphisms $\Diff (B,\omega)$ in order to translate
statements about maximal tori in $\Diff (B, \omega)$ into statements
about maximal tori in $(\Diff (P, \xi)$.  By a theorem of
Banyaga the group $\Diff (P, \alpha)$ of strict contactomorphism is a
central extension (possibly nontrivial) of $\Diff (B, \omega)$ by
$S^1$. However the group $\Diff (P, \xi )$ is much bigger than $\Diff
(P, \alpha)$.

\begin{remark}
When one talks about ``maximal tori'' in $\Diff(P, \xi)$ or in $\Diff
(B, \omega)$, one can mean two different things:
\begin{description}
\item[(1)] A torus $T$ in a group $G$ is {\em maximal} if for 
any torus $T' \subset G$ with $T\subset T'$ we have $T = T'$.

\item[(2a)]  A torus  $T \subset \Diff (B, \omega)$ satisfies 
$\dim T \leq \frac{1}{2} \dim B$ so a torus $T\subset \Diff (B,
\omega)$ of dimension $\frac{1}{2} \dim B$ is maximal.

\item[(2b)] Similarly, a torus  $T \subset \Diff (P, \xi)$ satisfies 
$\dim T \leq \frac{1}{2}( \dim P + 1)$  so a torus $T\subset \Diff (P,
\xi)$ of dimension $\frac{1}{2} (\dim P + 1)$ is maximal.
\end{description}

Since there exist compact simply connected symplectic four-manifolds
$(B, \omega)$ that admit no symplectic circle actions \cite{HK} (1)
and (2a) are quite different.  We will see that (1) and (2b) are
different as well, so we will use (1) as our definition of a maximal
torus.

\end{remark}

\noindent
The main result of the note is the following observation:

\begin{theorem}\label{thm1}
Let $(B, \omega)$ be a compact simply connected integral symplectic
manifold, $S^1 \to P \stackrel{\pi}{\to} B$ the principal circle
bundle with Euler class $[\omega]$, $\alpha \in\Omega^1 (P)$ a
connection 1-form with $d\alpha = \pi^* \omega$ and $\xi = \ker
\alpha$ the corresponding contact structure on $P$.  If $T \subset 
\Diff (B, \omega)$ is a maximal torus then its preimage $T'$ under the 
surjection $\tau : \Diff (P, \alpha) \to \Diff (B, \omega)$ is a maximal
torus in $\Diff (P, \xi)$.
\end{theorem}

\begin{remark}
The fact that the map $\tau :\Diff (P, \alpha) \to \Diff (B, \omega)$
in Theorem~\ref{thm1} above exists and is a surjective homomorphism is a
theorem of Banyaga \cite[Theorem~1]{B}.
\end{remark}
Our proof of Theorem~\ref{thm1} is a combination of two Lemmas
below. Note that there is a distinct $S^1$ in $\Diff (P, \alpha)
\subset \Diff (P, \xi)
\subset \Diff (P)$.  It is the circle action that makes $P$ a
principal $S^1$ bundle over $B$.  We will refer to this subgroup as
{\em the} $S^1$ in $\Diff (P)$.  
\begin{lemma}\label{lemma1}
The preimage of a torus $T\subset \Diff (B, \omega)$ in $\Diff (P,
\alpha)$ under $\tau$ is indeed a torus.  In fact the action of $T$ on
$B$ lifts to an action of $T$ on $P$ preserving $\alpha$ and commuting
with the action of the $S^1$. 
\end{lemma}
\begin{lemma}\label{lemma2}
If $H\subset \Diff (P, \xi)$ is a torus containing the $S^1$ then $H
\subset \Diff (P, \alpha)$.
\end{lemma}

Lemma~\ref{lemma1} is almost certainly not new and must exist
somewhere in the literature on pre-quantization of group actions.
Unfortunately I have been unable to find a good reference for the
result.  It is also possible that one can deduce it from
\cite[Theorem~1]{B}, but I am not whether  the group $\tau\inv (T)$ 
inherits from $\Diff (P, \alpha)$ the structure of a Lie group making
$\tau: \tau \inv (T) \to T $ a surjective Lie group homomorphism
(Banyaga is not very explicit about the topologies and smooth
structures of the groups involved in \cite[Theorem~1]{B}).  It is easy
to prove that if a Lie group $G$ is a central extension of a torus $T$
by $S^1$ as a Lie group and not just as an abstract group, then $G$ is
$S^1 \times T$.  On the other hand, Lemma~\ref{lemma2} has a very easy
proof.

\begin{proof}[Proof of Lemma~\ref{lemma2}] Let $R$ denote the vector field 
generating the $S^1$ action on $P$.  Then $\alpha (R) = 1$ and $L_R
\alpha =0$ since $\alpha$ is a connection.  Since $H$ contains the $S^1$, for 
any vector $X$ in the Lie algebra $\fh$ of $H$ the vector field $X_P$
induced by $X$ on $P$ commutes with $R$.  Also, since $H\subset \Diff
(P, \xi = \ker \alpha)$, $L_{X_P} \alpha = f \alpha$ for some function
$f \in C^\infty (P)$ that may depend on $X$.  We want to show that
$f\equiv 0$.  Now $$ 0 = L_{X_P} 1 = L_{X_P} (\iota (R)\alpha) = \iota
(L_{X_P} R) \alpha +
\iota (R)\left( L_{X_P} \alpha \right) = 0 + \iota (R) f\alpha = f.
$$
\end{proof}
Lemma~\ref{lemma1} is a consequence of the following slightly more
general proposition.  We drop the assumption that $B$ is simply
connected and replace it with the assumption that the action of $T$ is
Hamiltonian.  One may also drop the assumption that $B$ is compact,
but we won't do it. I am grateful to Anton Alekseev for telling me how
to prove the proposition.
\begin{proposition}
Let $(B, \omega)$ be a compact integral symplectic manifold.  Let $S^1
\to P \stackrel{\pi}{\to} B$ be the principal $S^1$ bundle with first
Chern (Euler) class $c_1 (P) = [\omega]$.  Let $\alpha$ be the
connection 1-form on $P$ with $d\alpha = \pi^* \omega$.  Suppose there
is a Hamiltonian action of a torus $T$ on $B$ with an associated
moment map $\Phi : B \to \ft^*$.

Then there is an action of $T$ on $P$ preserving $\alpha$ and making
$\pi$ equivariant.
\end{proposition}

\begin{proof}
It is well known how to lift the action of the Lie algebra $\ft$ of
$T$ on $B$ to an action on $P$ preserving $\alpha$ (see, for example
\cite{Sn}):  Given $X\in \ft$ the induced vector field $X_P$ on $P$ is
defined by 
\begin{equation}\label{eq}
X_P : = X_B ^h - (\pi^* \Phi^X) R, 
\end{equation}
where $X_B$ is the vector field induced by $X$ on $B$, $X_B^h$ denotes
its horizontal lift to $P$, $\Phi^X = \langle \Phi, X\rangle$ is the
$X$-component of the moment map $\Phi$, and $R$ is the vector field
generating the action of $S^1 = \R/\Z$.  Thus we have an action of the
universal cover $\tilde {T}$ of $T$ on $P$.  The point of the
proposition is that for a suitable normalization of the moment map
$\Phi$, the action of $\tilde{T}$ descends to an action of $T$.  Note
that the proposition is false if the group in question is {\em not} a
torus: there is no way to lift the standard action of $SO(3)$ on $S^2$
to an action on $S^3$.

Our normalization is as follows.  Since $B$ is compact and the action
of $T$ is Hamiltonian, the set of $T$-fixed points $B^T$ is non-empty.
Pick a  point $b_0 \in B^T$ and normalize $\Phi$ by requiring
that $\Phi (b_0) = 0$.  We claim for any vector $X$ in the integral
lattice $\Z_T :=\ker \{\exp: \ft \to T\}$ the flow $\{\exp t X_P \}$
defined by (\ref{eq}) is periodic of period 1.  Indeed, since the
vector fields $X_B^h$ and $(\pi^* \Phi^X) R$ commute, $$
\exp t X_P = (\exp t X_B ^h)\circ \exp ( -t (\pi^* \Phi^X) R).
$$
For a point $p\in P$,
$$
\exp ( -t (\pi^* \Phi^X) R) (p) = e^{- 2\pi i t \Phi^X (b)} ,
$$ 
where $b = \pi (p)$ and where we identified $\R/\Z $ with $U(1)$ by
$\theta \!\mod \Z \mapsto e^{2\pi i \theta}$.

The curve $\gamma (t) = (\exp t X_B) (b)$ is a loop in $B$ since $X\in
\Z_T$.  Hence $(\exp t X_B ^h) (p) = H(\gamma)\cdot p$ where
$H(\gamma)$ denotes the holonomy of $\gamma$.  On the other hand, 
$$
H(\gamma) = e^{2\pi i \int _D \omega} 
$$ 
for any disk $D \subset B$ with boundary $\gamma$. (Note that if $D'$
is another disk with boundary $\gamma$, then $\int_D \omega - \int
_{D'} \omega = \int _{S^2} \omega \in \Z$, since $[\omega]$ is
integral, and consequently $e^{2\pi i \int _D \omega}$ is
well-defined.)\,

The curve $\gamma$ always bounds a disk: since $B$ is connected, there
is a path $\tau : [0, 1] \to B$ with $\tau (0) = b_0 $ and $\tau (1) =
b$.  The disk
$$
D_\tau  
= \{ (\exp tX_B)\cdot \tau (s) \mid 0\leq t \leq 1, \,\, 0 \leq s \leq 1\}
$$
is a desired disk.  Moreover,
$$
\int _{D_\tau} \omega = - \int _{\tau ([0, 1])} (\iota (X_B) \omega)=
\int _0^1 d\Phi^X (\tau (s))\,ds = \Phi^X (\tau (1)) - \Phi ^X (\tau (0))
= \Phi^X (b). 
$$
Thus $(\exp X_B ^h) \cdot p = H (\gamma)\cdot p = e^{2\pi i \Phi ^X
(b)} \cdot p$, and therefore 
$$ 
(\exp X_P)(p) = \exp X_B^h \cdot
e^{-2\pi i \Phi ^X (b)} \cdot p = e^{2\pi i \Phi ^X (b)}e^{-2\pi i \Phi
^X (b)} \cdot p = p.  
$$
\end{proof}

\section{Examples}

\begin{example}
Hausmann and Knutson \cite{HK} constructed a symplectic form $\omega$
on $B = \C P \# 3 \overline{\C P^2}$ which admits no Hamiltonian
(hence symplectic) circle actions.  The form $\omega$ may be taken to
be integral.  The manifold $(B, \omega)$ is a pentagon space.
Consider the corresponding contact manifold $(P, \xi = \ker \alpha)$,
where as above $\pi: P\to B$ is the principal $S^1$ bundle with $c_1
(P) = [\omega]$ and $\alpha$ is a connection 1-form with $d\alpha =
\pi^* \omega$.  By Theorem~\ref{thm1} the $S^1$ in the contactomorphism 
group $\Diff (P, \xi)$ is a maximal torus.  Note that $\dim P =5$, so
the maximal possible dimension of a maximal torus in $\Diff (P, \xi)$
is 3.
\end{example}

\begin{example}
Hausmann and Tolman \cite{HT} constructed a number of polygon spaces
$(B, \omega)$ with the property that the symplectomorphism group
$\Diff (B, \omega)$ contains maximal tori of different dimensions.
For example the group of symplectomorphisms of the heptagon space
$\Pol (1, 1, 2, 2, 3, 3, 3)$ (we use the notation of \cite{HT})
contains maximal tori of dimensions 2, 3 and 4.  Hence the
contactomorphism group $\Diff (P, \xi)$ of the corresponding principal
circle bundle $P\to \Pol (1, 1, 2, 2, 3, 3, 3)$ contains maximal tori
of dimension $2+1$, $3+1$ and $4+1$.
\end{example}

\section*{Acknowledgments}
The work on this note was supported by the Swiss NSF at the
University of Geneva in May 2002.  I am grateful to the University for
its hospitality.  I thank J.-C. Hausmann and A. Alekseev for their
help.

\end{document}